\documentclass[12pt]{amsart}
\usepackage{amsmath, amssymb}
\usepackage{amsfonts}
\usepackage{amscd}
\usepackage[arrow,matrix,curve,cmtip,ps]{xy}

\usepackage{amsthm}

\allowdisplaybreaks

\newtheorem{theorem}{Theorem}[section]
\newtheorem{lemma}[theorem]{Lemma}
\newtheorem{proposition}[theorem]{Proposition}

\newtheorem{corollary}[theorem]{Corollary}

\newtheorem*{theorem*}{Theorem}
\theoremstyle{definition}
\newtheorem{remark}[theorem]{Remark}
\newtheorem{definition}[theorem]{Definition}
\newtheorem{example}[theorem]{Example}

\numberwithin{equation}{section}

\newcommand{\Z}{\mathbb{Z}}
\newcommand{\N}{\mathbb{N}}

\newcommand{\C}{\mathbb{C}}
\newcommand{\T}{\mathbb{T}}
\newcommand{\G}{\mathcal{G}}
\newcommand{\cH}{\mathcal{H}}
\newcommand{\A}{\mathfrak{A}}
\newcommand{\B}{\mathfrak{B}}

\newcommand{\cC}{\mathcal{C}}
\newcommand{\cO}{{\mathcal O}}
\newcommand{\cQ}{{\mathcal Q}}

\newcommand{\Pow}{\mathcal{P}}
\newcommand{\tDelta}{\widetilde{\Delta}}
\newcommand{\tOmega}{\widetilde{\Omega}}
\newcommand{\rg}{\textnormal{rg}}
\newcommand{\coker}{\operatorname{coker}}

\DeclareMathOperator{\cspa}{\overline{span}}
\newcommand{\rom}{\renewcommand{\labelenumi}{{\rm (\roman{enumi})}}%
\renewcommand{\itemsep}{0pt}}
\newcommand{\Ca}{$C^*$-al\-ge\-bra }
\newcommand{\CA}{$C^*$-al\-ge\-bra}
\newcommand{\Csa}{$C^*$-sub\-al\-ge\-bra }
\newcommand{\CsA}{$C^*$-sub\-al\-ge\-bra}
\newcommand{\shom}{$*$-ho\-mo\-mor\-phism }
\newcommand{\shoM}{$*$-ho\-mo\-mor\-phism}
\begin{document}
\title[Ultragraphs and quivers]{Ultragraph $C^*$-algebras via topological quivers}

\author[Katsura]{Takeshi Katsura}
\address{T. Katsura, Department of Mathematics\\ Keio University\\
Yokohama, 223-8522\\ JAPAN}
\email{katsura@math.keio.ac.jp}

\author[Muhly]{Paul S. Muhly}
\address{P. Muhly, Department of Mathematics\\ University of Iowa\\
Iowa City\\ IA 52242-1419\\ USA} \email{pmuhly@math.uiowa.edu}

\author[Sims]{Aidan Sims}
\address{A. Sims, School of Mathematics and Applied Statistics\\
University of Wollongong\\
NSW 2522\\
AUSTRALIA} \email{asims@uow.edu.au}

\author[Tomforde]{Mark Tomforde}
\address{M. Tomforde \\ Department of Mathematics\\ University of Houston\\
Houston \\ TX 77204-3008\\ USA} \email{tomforde@math.uh.edu}

\thanks{The first author was supported by JSPS, the second author was supported by NSF Grant DMS-0070405, the third author was supported by the Australian Research Council, and the fourth author was supported by NSF Postdoctoral Fellowship DMS-0201960.}

\date{\today}
\subjclass[2000]{Primary 46L05}

\keywords{$C^*$-algebras, ultragraphs, topological graphs}

\begin{abstract}
Given an ultragraph in the sense of Tomforde, we construct a
topological quiver in the sense of Muhly and Tomforde in such a way
that the universal $C^*$-algebras associated to the two objects
coincide. We apply results of Muhly and Tomforde for topological
quiver algebras and of Katsura for topological graph $C^*$-algebras
to study the $K$-theory and gauge-invariant ideal structure of
ultragraph $C^*$-algebras.
\end{abstract}

\maketitle

\section{Introduction}

Our objective in this paper is to show how the theory of
\emph{ultragraph \CA s}, first proposed by Tomforde in
\cite{Tom, Tom2}, can be formulated in the context of topological
graphs \cite{Kat} and topological quivers \cite{MT} in a fashion that
reveals the $K$-theory and ideal theory (for gauge-invariant ideals)
of these algebras. The class of graph \CA s has attracted
enormous attention in recent years. The graph \Ca
associated to a directed graph $E$ is generated by projections $p_v$
associated to the vertices $v$ of $E$ and partial isometries $s_e$
associated to the edges $e$ of $E$. Graph \CA s, which, in turn,
are a generalization of the Cuntz-Krieger algebras of \cite{CK}, were
first studied using groupoid methods \cite{KPRR, KPR}. An artifact of
the initial groupoid approach is that the original theory was
restricted to graphs which are \emph{row-finite} and have \emph{no
sinks} in the sense that each vertex emits at least one and at most
finitely many edges\footnote{A refinement of the analysis in
\cite{KPRR}, due to Paterson \cite{aP02}, extends the groupoid
approach to cover non-row-finite graphs. One can find a groupoid
approach that handles sinks and extends the whole theory to higher
rank graphs in \cite{FMY05}. A groupoid approach to ultragraph
\CA s may be found in \cite{MMp06}.}.

The connection between Cuntz-Krieger algebras and graph
\CA s is that each directed graph can be described in terms
of its edge matrix, which is a $\{0,1\}$-matrix indexed by the edges
of the graph; a $1$ in the $(e,f)$ entry indicates that the range of
$e$ is equal to the source of $f$. The terminology \emph{row-finite}
refers to the fact that in any row of the edge matrix of a row-finite
graph, there are at most finitely many nonzero entries.

The two points of view, graph and matrix, led to two versions of
Cuntz-Krieger theory for non-row-finite objects. In \cite{FLR},
\CA s were associated to arbitrary graphs in such a way that
the construction agrees with the original definition in the
row-finite case. In \cite{EL} \CA s --- now called Exel-Laca
algebras --- were associated to arbitrary $\{0,1\}$ matrices, once
again in such a way that the definitions coincide for row-finite
matrices. The fundamental difference between the two classes of
algebras is that a graph \Ca is generated by a collection
containing a partial isometry for each edge and a projection for each
vertex, while an Exel-Laca algebra is generated by a collection
containing a partial isometry for each row in the matrix (and in the
non-row-finite case there are rows in the matrix corresponding to an
infinite collection of edges with the same source vertex).  Thus,
although these two constructions agree in the row-finite case, there
are \CA s of non-row-finite graphs that are not isomorphic
to any Exel-Laca algebra, and there are Exel-Laca algebras of
non-row-finite matrices that are not isomorphic to the \Ca
of any graph \cite{Tom2}.

In order to bring graph \CA s of non-row-finite graphs and
Exel-Laca algebras together under one theory, Tomforde introduced the
notion of an ultragraph and described how to associate a
\Ca to such an object \cite{Tom, Tom2}. His analysis not
only brought the two classes of \CA s under one rubric, but
it also showed that there are ultragraph \CA s that belong
to neither of these classes.  Ultragraphs are basically directed
graphs in which the range of each edge is non-empty \emph{set} of
vertices rather than a single vertex --- thus in an ultragraph each
edge points from a single vertex to a set of vertices, and directed
graphs are the special case where the range of each edge is a
singleton set.  Many of the fundamental results for graph
\CA s, such as the well-known Cuntz-Krieger Uniqueness
Theorem and the Gauge-Invariant Uniqueness Theorem, can be proven in
the setting of ultragraphs \cite{Tom}. However, other results, such
as $K$-theory computations and ideal structure are less obviously
amenable to traditional graph \Ca techniques.

Recently, Katsura \cite{Kat} and Muhly and Tomforde \cite{MT} studied
the notions of topological graphs and topological quivers,
respectively. These structures consist of second countable locally
compact Hausdorff spaces $E^0$ and $E^1$ of vertices and edges
respectively with range and source maps $r,s : E^1 \to E^0$ which
satisfy appropriate topological hypotheses. The main point of
difference between the two (apart from a difference in edge-direction
conventions) is that in a topological graph the  source map is
assumed to be a local homeomorphism so that $s^{-1}(v)$ is discrete,
whereas in a topological quiver the range map (remember the
edge-reversal!) is only assumed to be continuous and open and a
system $\lambda = \{ \lambda_v \}_{v \in E^0}$ of Radon measures
$\lambda_v$ on $r^{-1}(v)$ satisfying some natural conditions (see
\cite[Definition~3.1]{MT}) is supplied as part of the data. It is
worth pointing out that given $E^0, E^1, r$ and $s$, with $r$ open,
such a system of Radon measures will always exist. A topological
graph can be regarded as a topological quiver by reversing the edges
and taking each $\lambda_v$ to be counting measure; the topological
graph \Ca and the topological quiver \Ca then
coincide. One can regard an ordinary directed graph as either a
topological graph or a topological quiver by endowing the edge and
vertex sets with the discrete topology, and  then the topological
graph \Ca and topological quiver algebra coincide with the
original graph \CA.

In this article we show how to build a topological quiver $\cQ(\G)$
from an ultragraph $\G$ in such a way that the ultragraph
\Ca $C^*(\G)$ and the topological quiver algebra
$C^*(\cQ(\G))$ coincide. We then use results of \cite{Kat} and
\cite{MT} to compute the $K$-theory of $C^*(\G)$, to produce a
listing of its gauge-invariant ideals, and to provide a version of
Condition~(K) under which all ideals of $C^*(\G)$ are
gauge-invariant.

It should be stressed that the range map in $\cQ(\G)$ is always a
local homeomorphism, so $\cQ(\G)$ can equally be regarded as a
topological graph; indeed our analysis in some instances requires
results regarding topological graphs from \cite{Kat} that have not
yet been generalized to topological quivers. We use the notation and
conventions of topological quivers because the edge-direction
convention for quivers in \cite{MT} is compatible with the
edge-direction convention for ultragraphs \cite{Tom, Tom2}.

\vskip0.5em

The paper is organized as follows: In Section~\ref{AG-sec} we
describe the commutative \Ca $\A_\G \subset C^*(\G)$
generated by the projections $\{p_A : A \in \G^0\}$. In
Section~\ref{C*alg-sec} we provide two alternative formulations of
the defining relations among the generators of an ultragraph
\Ca which will prove more natural in our later analysis. In
Section~\ref{OmegaG-sec} we describe the spectrum of $\A_\G$. We use
this description in Section~\ref{TopQu&K-sec} to define the quiver
$\cQ(\G)$, show that its \Ca is isomorphic to $C^*(\G)$,
and compute its $K$-theory in terms of the structure of $\G$ using
results from \cite{Kat}. In Section~\ref{Ideal-sec} we use the
results of \cite{MT} to produce a listing of the gauge-invariant
ideals of $C^*(\G)$ in terms of the structure of $\G$, and in
Section~\ref{CondK-sec} we use a theorem of \cite{Kat2} to
provide a condition on $\G$ under which all ideals of $C^*(\G)$ are
gauge-invariant.

\section{The commutative $C^*$-algebra $\A_\G$ and its representations} \label{AG-sec}

For a set $X$, let $\Pow(X)$ denote the collection of all subsets of
$X$.

\begin{definition}
An \emph{ultragraph} $\G = (G^0, \G^1, r,s)$ consists of a countable
set of vertices $G^0$, a countable set of edges $\G^1$, and functions
$s\colon \G^1\rightarrow G^0$ and $r\colon \G^1 \rightarrow
\Pow(G^0)\setminus\{\emptyset\}$.
\end{definition}

\begin{definition}
For an ultragraph $\G=(G^0, \G^1, r,s)$, we denote by $\A_\G$ the
\Csa of $\ell^\infty(G^0)$ generated by the point masses
$\{\delta_{v}:v\in G^0\}$ and the characteristic functions
$\{\chi_{r(e)}: e\in\G^1\}$.
\end{definition}

Let us fix an ultragraph $\G =(G^0, \G^1, r,s)$, and consider the
representations of $\A_\G$.

\begin{definition}
For a set $X$, a subcollection $\cC\/$ of $\Pow(X)$ is called a {\em
lattice} if
\begin{enumerate}
\item[(i)] $\emptyset\in\cC$
\item[(ii)] $A\cap B\in\cC$ and $A\cup B\in\cC$ for all $A,B\in\cC$.
\end{enumerate}
An {\em algebra} is a lattice $\cC$ that also satisfies the
additional condition
\begin{enumerate}
\item[(iii)] $A\setminus B\in \cC$ for all $A,B\in\cC$.
\end{enumerate}
\end{definition}

\begin{definition}
For an ultragraph $\G = (G^0, \G^1, r,s)$, we let $\G^0$ denote the
smallest algebra in $\Pow(G^0)$ containing the singleton sets and the
sets  $\{r(e) : e \in \G^1\}$.
\end{definition}

\begin{remark}
In \cite{Tom}, $\G^0$ was defined to be the smallest lattice --- not
algebra --- containing the singleton sets and the sets  $\{r(e) : e
\in \G^1\}$. The change to the above definition causes no problem
when defining Cuntz-Krieger $\G$-families (see the final paragraph of
Section \ref{C*alg-sec}).  Furthermore, this new definition is convenient for us in a variety of situations:  It relates $\G^0$ to the $C^*$-algebra $\A_\G$ described in Proposition~\ref{prop:A_G}, it allows us too see immediately that the set $r(\lambda, \mu)$ of Definition~\ref{def: r-lam-mu-def} is in $\G^0$, and --- most importantly --- it aids in our description of the gauge-invariant ideals in Definition~\ref{def:ideal with new def} and Lemma~\ref{lem:ideal corresp}.  For additional justification for the change in definition, we refer the reader to \cite[Section~2]{MMp06}.
\end{remark}

\begin{proposition}\label{prop:A_G}
We have $\G^0=\{A\subset G^0 : \chi_A\in\A_\G\}$ and
\begin{equation}\label{eq:span expr}
\A_\G=\cspa\{\chi_A : A\in \G^0\}.
\end{equation}
\end{proposition}

\begin{proof}
We begin by proving~\eqref{eq:span expr}. Since $\{A\subset G^0 :
\chi_A\in\A_\G\}$ is an algebra containing $\{v\}$ and $r(e)$, we
have $\G^0\subset \{A\subset G^0 : \chi_A\in\A_\G\}$. Hence we get
$$\A_\G\supset\cspa\{\chi_A : A\in \G^0\}.$$
Since $\G^0$ is closed under intersections, the set $\cspa\{\chi_A :
A\in \G^0\}$ is closed under multiplication, and hence is a \Ca
containing $\{\delta_v\}$ and $\{\chi_{r(e)}\}$. Hence
$\A_\G\subset\cspa\{\chi_A : A\in \G^0\}$,
establishing~\eqref{eq:span expr}.

We must now show that $\G^0=\{A\subset G^0 : \chi_A\in\A_\G\}$. We
have already seen $\G^0 \subset \{A\subset G^0 : \chi_A\in\A_\G\}$.
Let $A\subset G^0$ with $\chi_A\in\A_\G$. By~\eqref{eq:span expr}, we
have $\|\chi_A - \sum_{k=1}^m z_k\chi_{A_k}\| < 1/2$ for some
$A_1,\ldots,A_m \in \G^0$ and $z_1,\ldots, z_m \in \C$; moreover,
since $\G^0$ is an algebra, we may assume that $j \not= k$ implies
$A_j \cap A_k = \emptyset$. But then $x \in A$ if and only if $x \in
A_k$ for some (unique) $k$ with $|1 - z_k| < 1/2$. That is,
$A=\bigcup_{|1 - z_k|<1/2}A_k\in \G^0$.
\end{proof}

\begin{definition}
A \emph{representation} of a lattice $\cC$ on a \Ca $\B$ is a
collection of projections $\{p_A\}_{A\in \cC}$ of $\B$ satisfying
$p_\emptyset = 0$, $p_A p_B = p_{A \cap B}$, and $p_{A\cup B} = p_A +
p_B - p_{A \cap B}$ for all $A,B \in \cC$.
\end{definition}

When $\cC$ is an algebra, the last condition of a representation can
be replaced by the equivalent condition that $p_{A\cup B} = p_A +
p_B$ for all $A,B \in \cC$ with $A \cap B=\emptyset$.

Note that we define representations of lattices here, rather than just of
algebras, so that our definition of $C^*(\G)$ agrees with the original
definition given in \cite{Tom} (see the final paragraph of
Section~\ref{C*alg-sec}).

\begin{definition} \label{def: r-lam-mu-def}
For a fixed ultragraph $\G = (G^0, \G^1, r, s)$ we define
\[
X=\big\{(\lambda,\mu): \text{$\lambda,\mu$ are finite subsets of
$\G^1$ with $\lambda\cap \mu=\emptyset$ and
$\lambda\neq\emptyset$}\big\}.
\]
For $(\lambda,\mu)\in X$, we define $r(\lambda,\mu)\subset G^0$ by
\[
r(\lambda,\mu):= \bigcap_{e\in \lambda}r(e)\setminus \bigcup_{f\in
\mu}r(f).
\]
\end{definition}

\begin{definition}\label{dfn:Cond(EL)}
Let $\G =(G^0, \G^1, r,s)$ be an ultragraph. A collection of
projections $\{p_v\}_{v\in G^0}$ and $\{q_e\}_{e \in \G^1}$ is said
to satisfy {\em Condition (EL)} if the following hold:
\begin{enumerate}
\item the elements of $\{p_v\}_{v\in G^0}$ are pairwise orthogonal,
\item the elements of $\{q_e\}_{e \in \G^1}$ pairwise commute,
\item $p_v q_e= p_v$ if $v\in r(e)$, and $p_v q_e=0$ if $v \notin r(e)$,
\item $\prod_{e\in \lambda}q_e\prod_{f\in \mu}(1-q_f)
=\sum_{v\in r(\lambda,\mu)}p_v$ for all $(\lambda,\mu)\in X$ with
$|r(\lambda,\mu)|<\infty$.
\end{enumerate}
\end{definition}

From a representation of $\G^0$, we get a collection satisfying
Condition~(EL). We prove a slightly stronger statement.

\begin{lemma}\label{lem:G'->(*)}
Let $\cC$ be a lattice in $\Pow(G^0)$ which contains the singleton
sets and the sets  $\{r(e) : e \in \G^1\}$, and let $\{p_A\}_{A \in
\cC}$ be a representation of $\cC$. Then the collection
$\{p_{\{v\}}\}_{v\in G^0}$ and $\{p_{r(e)}\}_{e \in \G^1}$ satisfies
Condition~(EL).
\end{lemma}

\begin{proof}
From the condition $p_\emptyset = 0$ and $p_A p_B = p_{A \cap B}$, it
is easy to show that the collection satisfies the conditions (1), (2)
and (3) in Definition~\ref{dfn:Cond(EL)}. To see condition~(4) let
$(\lambda,\mu)\in X$ with $|r(\lambda,\mu)|<\infty$. Define
$A,B\subset G^0$ by $A=\bigcap_{e\in \lambda}r(e)$ and
$B=\bigcup_{f\in \mu}r(f)$. Then we have $A,B\in \cC$, and from the
definition of a representation, we obtain
\[
p_A=\prod_{e\in \lambda}p_{r(e)} \qquad \text{ and } \qquad
1-p_B=\prod_{f\in \mu}(1-p_{r(f)}).
\]
Since $r(\lambda,\mu)$ is a finite set, $r(\lambda,\mu)\in \cC$ and
$p_{r(\lambda,\mu)}=\sum_{v\in r(\lambda,\mu)}p_{\{v\}}$. Also,
because $r(\lambda,\mu)=A\setminus B$, we obtain $r(\lambda,\mu)\cup
B=A\cup B$ and $r(\lambda,\mu)\cap B=\emptyset$. Hence $p_{A\cup
B}=p_{r(\lambda,\mu)}+p_B$. Since $p_{A\cup B}=p_A+p_B-p_{A\cap B}$,
we get $p_{r(\lambda,\mu)}=p_A-p_{A\cap B}$. Hence we have,
\[
\sum_{v\in r(\lambda,\mu)}p_{\{v\}}=p_A-p_{A\cap B}=p_A(1-p_B)
=\prod_{e\in \lambda}p_{r(e)}\prod_{f\in \mu}(1-p_{r(f)}).
\]
\end{proof}

We will prove that from a collection satisfying Condition~(EL), we
can construct a \shom from $\A_\G$ onto the $C^*$-subalgebra
generated by that collection. To this end, we fix a listing $\G^1=
\{e_i\}^\infty_{i=1}$, and for each positive integer $n$ define a
\Csa $\A_\G^{(n)}$ of $\A_\G$ to be the \Ca generated by
$\{\delta_{v}:v\in G^0\}$ and $\{\chi_{r(e_i)}:i=1,2,\ldots,n\}$.
Note that the union of the increasing family
$\{\A_\G^{(n)}\}_{n=1}^\infty$ is dense in $\A_\G$.

\begin{definition}
Let $n$ be a positive integer. Let $0^n := (0,0,\ldots,0)\in
\{0,1\}^{n}$. For $\omega=(\omega_1,\omega_2,\ldots,\omega_n) \in
\{0,1\}^n\setminus \{0^n\}$, we set
\[
r(\omega):=\bigcap_{\omega_i=1}r(e_i)\setminus
\bigcup_{\omega_j=0}r(e_j).
\]
\end{definition}

\begin{lemma}
Let $n$ be a positive integer. For each $\omega\in \{0,1\}^{n}$, we
define $\lambda_\omega,\mu_\omega\subset\G^1$ by
$\lambda_\omega=\{e_i:\omega_i=1\}$ and
$\mu_\omega=\{e_i:\omega_i=0\}$. Then the map
\[
\omega\mapsto (\lambda_\omega,\mu_\omega)
\]
is a bijection between $\{0,1\}^n\setminus \{0^n\}$ and
$\big\{(\lambda,\mu)\in X: \lambda\cup\mu=\{e_1,\ldots,e_n\}\big\}$,
and we have $r(\omega)=r(\lambda_\omega,\mu_\omega)$.
\end{lemma}
\begin{proof}
The map $\omega \mapsto (\lambda_\omega,\mu_\omega)$ is a bijection
because $(\lambda,\mu) \mapsto \chi_\lambda$ provides an inverse, and
$r(\omega)=r(\lambda_\omega,\mu_\omega)$ by definition.
\end{proof}

\begin{definition}
We define $\Delta_n :=\big\{\omega\in \{0,1\}^n\setminus \{0^n\}:
|r(\omega)|=\infty\big\}.$
\end{definition}

\begin{lemma}\label{lem:r(omega)}
For each $i=1,2,\ldots,n$, the set $r(e_i)$ is a disjoint union of
the infinite sets $\{r(\omega)\}_{\omega \in \Delta_n,\omega_i=1}$
and the finite set $\bigcup_{\omega \notin
\Delta_n,\omega_i=1}r(\omega)$.
\end{lemma}
\begin{proof}
First note that $r(\omega) \cap r(\omega') = \emptyset$ for distinct
$\omega,\omega' \in \{0,1\}^{n}\setminus\{0^n\}$. For $v \in r(e_i)$,
define $\omega^v \in \{0,1\}^n$ by $\omega^v_j = \chi_{r(e_j)}(v)$
for $1 \le j \le n$. Since $v \in r(e_i)$, $\omega^v \not= 0^n$, and
$v \in r(\omega^v)$ by definition. Hence
\[
r(e_i) = \bigcup_{v \in r(e_i)} r(\omega^v) = \bigcup_{\omega_i = 1}
r(\omega)
\]
Since $r(\omega)$ is a finite set for $\omega\in \{0,1\}^n\setminus
\{0^n\}$ with $\omega\notin\Delta_n$, the result follows.
\end{proof}

\begin{proposition} \label{prop:A-G-n generated}
The \Ca $\A_\G^{(n)}$ is generated by $\{\delta_{v}:v\in G^0\} \cup
\{\chi_{r(\omega)}:\omega\in \Delta_n\}$.
\end{proposition}
\begin{proof}
For each $\omega\in \Delta_n$, we have
\[
\chi_{r(\omega)}=\prod_{\omega_i=1}\chi_{r(e_i)}
\prod_{\omega_j=0}\big(1-\chi_{r(e_j)}\big)\in \A_\G^{(n)},
\]
giving inclusion in one direction.  It follows from Lemma~\ref{lem:r(omega)}
that the generators of $\A_\G^{(n)}$ all belong to the \Ca generated
by $\{\delta_{v}:v\in G^0\} \cup \{\chi_{r(\omega)}:\omega\in
\Delta_n\}$, establishing the reverse inclusion.
\end{proof}

For each $\omega\in \Delta_n$, the \Csa of $\A_\G^{(n)}$ generated by
$\{\delta_v:v\in r(\omega)\}$ and $\chi_{r(\omega)}$ is isomorphic to
the unitization of $c_0(r(\omega))$. Since the \Ca $\A_\G^{(n)}$ is a
direct sum of such \CsA s indexed by the set $\Delta_n$ and the \Csa
$c_0\big(G^0\setminus\bigcup_{\omega\in \Delta_n}r(\omega)\big)$
(recall that the $r(\omega)$'s are pairwise disjoint), we have the
following:

\begin{lemma}\label{lem:An->B}
For two families $\{p_v\}_{v\in G^0}$ and $\{q_\omega\}_{\omega\in
\Delta_n}$ of mutually orthogonal projections in a $C^*$-algebra $\B$ satisfying
\[
p_vq_\omega=\begin{cases}p_v& \text{if $v\in r(\omega)$,}\\
0& \text{if $v\notin r(\omega)$,}\end{cases}
\]
there exists a $*$-homomorphism $\pi_n\colon \A_\G^{(n)}\to \B$ with
$\pi_n(\delta_v)=p_v$ for $v\in G^0$ and
$\pi_n(\chi_{r(\omega)})=q_\omega$ for $\omega\in \Delta_n$.
\end{lemma}

\begin{proposition}\label{prop:rep<->EL}
Let $\G=(G^0, \G^1, r,s)$ be an ultragraph, and $\B$ be a \CA.
Then there exist natural bijections among the following sets:
\begin{enumerate}
\rom
\item the set of \shoM s from $\A_\G$ to $\B$,
\item the set of representations of $\G^0$ on $\B$,
\item the set of collections of projections in $\B$ satisfying Condition~(EL).
\end{enumerate}
Specifically, if $\pi\colon \A_\G \to \B$ is a \shoM, then $p_A :=
\pi(\chi_A)$ for $A \in \G^0$ gives a representation of $\G^0$; if
$\{p_A\}_{A \in \G^0}$ is a representation of $\G^0$ on $\B$, then
$\{p_{\{v\}}\}_{v\in G^0} \cup \{p_{r(e)}\}_{e \in \G^1}$ satisfies
Condition~(EL); and if a collection of projections $\{p_v\}_{v\in
G^0} \cup \{q_e\}_{e \in \G^1}$ satisfies Condition~(EL), then there
exists a unique \shoM\ $\pi\colon \A_\G \to \B$ such that
$\pi(\delta_{v})=p_v$ and $\pi(\chi_{r(e)})=q_e$ for all $v\in G^0$
and $e \in \G^1$.
\end{proposition}
\begin{proof}
Clearly we have the map from (i) to (ii), and by Lemma
\ref{lem:G'->(*)} we have the map from (ii) to (iii). Suppose that
$\{p_v\}_{v\in G^0}$ and $\{q_e\}_{e \in \G^1}$ is a collection of
projections satisfying Condition (EL). Fix a positive integer $n$.
For each $\omega\in\{0,1\}^n\setminus\{0^n\}$, we define
$q_\omega=\prod_{\omega_i=1}q_{e_i}
\prod_{\omega_j=0}(1-q_{e_j})\in\B$. Then
$\{q_\omega:\omega\in\{0,1\}^n\setminus\{0^n\}\}$ is mutually
orthogonal. By Definition~\ref{dfn:Cond(EL)}(3), we have
\[
p_vq_\omega=\begin{cases}p_v& \text{if $v\in r(\omega)$,}\\
0& \text{if $v\notin r(\omega)$.}\end{cases}
\]
Hence by Lemma~\ref{lem:An->B}, there exists a \shom $\pi_n\colon
\A_\G^{(n)}\to \B$ such that $\pi_n(\delta_v)=p_v$ for $v\in G^0$ and
$\pi_n(\chi_{r(\omega)})=q_\omega$ for $\omega\in \Delta_n$. For
$\omega\in \{0,1\}^n\setminus\{0^n\}$ with $|r(\omega)|<\infty$, we
have $\pi_n(\chi_{r(\omega)})=\sum_{v\in r(\omega)}p_v=q_\omega$ by
Definition~\ref{dfn:Cond(EL)}(4). Hence we obtain
\[
\pi_n(\chi_{r(e_i)})=\pi_n\bigg(
\sum_{\substack{\omega \in \{0,1\}^n \notag \\
\omega_i=1}}\chi_{r(\omega)}\bigg) =\sum_{\substack{\omega \in
\{0,1\}^n \notag \\ \omega_i=1}}q_\omega =q_{e_i},
\]
for $i=1,\ldots,n$. Thus for each $n$, the \shom $\pi_n\colon
\A_\G^{(n)}\to \B$ satisfies $\pi_n(\delta_v)=p_v$ for $v\in G^0$ and
$\pi_n(\chi_{r(e_i)})=q_{e_i}$ for $i=1,\ldots,n$. Since there is at
most one $*$-homomorphism of $\A_\G^{(n)}\to \B$ with this property,
the restriction of the $*$-homomorphism $\pi_{n+1} : \A_\G^{(n+1)}
\to \B$ to $\A_\G^{(n)}$ coincides with $\pi_n$. Hence there is a
\shom $\pi\colon \A_\G\to \B$ such that $\pi(\delta_v)=p_v$ for $v\in
G^0$ and $\pi(\chi_{r(e)})=q_{e}$ for $e\in \G^1$. The \shom $\pi$ is
unique because $\A_\G$ is generated by $\{\delta_{v}:v\in G^0\} \cup
\{\chi_{r(e)}: e\in\G^1\}$.
\end{proof}

\begin{corollary}\label{cor:cond(4)}
Let $\G^1= \{e_i\}^\infty_{i=1}$ be some listing of the edges of
$\G$. To check that a family of projections $\{ p_v \}_{v \in G^0}
\cup \{q_e \}_{e \in \G^1}$ satisfies Condition~(EL), it suffices to
verify that Definition~\ref{dfn:Cond(EL)}(1)--(3) hold and that (4)
holds for $(\lambda,\mu)\in X$ with $|r(\lambda,\mu)|<\infty$ and
$\lambda\cup\mu=\{e_1,\ldots,e_n\}$ for some $n$.
\end{corollary}

We conclude this section by computing the $K$-groups of the \Ca
$\A_\G$.

\begin{definition}
For an ultragraph $\G = (G^0, \G^1, r,s)$, we denote by $Z_\G$ the
(algebraic) subalgebra of $\ell^\infty(G^0,\Z)$ generated by
$\{\delta_{v}:v\in G^0\} \cup \{\chi_{r(e)}: e\in\G^1\}$.
\end{definition}

An argument similar to the proof of Proposition~\ref{prop:A_G} shows
that
\[
Z_\G=\bigg\{\sum_{k=1}^n z_k \chi_{A_k} : n\in\N, z_k\in\Z, A_k\in
\G^0\bigg\}.
\]

\begin{proposition}\label{prop:K_0}
We have $K_0(\A_\G)\cong Z_\G$ and $K_1(\A_\G)=0$.
\end{proposition}
\begin{proof}
For each $n \in \N\setminus\{0\}$, let $Z_\G^{(n)}$ be the subalgebra
of $\ell^\infty(G^0,\Z)$ generated by $\{\delta_{v}:v\in G^0\} \cup
\{\chi_{r(e_i)}: i=1,2,\ldots,n\}$. By an argument similar to the
paragraph following Proposition~\ref{prop:A-G-n generated}, we see
that $Z_\G^{(n)}$ is a direct sum of the unitizations (as algebras)
of $c_0(r(\omega),\Z)$'s and $c_0\big(G^0\setminus\bigcup_{\omega\in
\Delta_n}r(\omega), \Z\big)$. Hence the description of $\A_\G^{(n)}$
in the paragraph following Proposition~\ref{prop:A-G-n generated}
shows that there exists an isomorphism $K_0(\A_\G^{(n)})\to
Z_\G^{(n)}$ which sends $[\delta_{v}],[\chi_{r(\omega)}]\in
K_0(\A_\G^{(n)})$ to $\delta_{v},\chi_{r(\omega)}\in Z_\G^{(n)}$. By
taking inductive limits, we get an isomorphism $K_0(\A_\G)\to Z_\G$
which sends $[\chi_A]\in K_0(\A_\G)$ to $\chi_{A}\in Z_\G$ for
$A\in\G^0$. That $K_1(\A_\G)=0$ follows from the fact that
$K_1(\A_\G^{(n)})=0$ for each $n$, and by taking direct limits.
\end{proof}

\begin{remark}
It is not difficult to see that the isomorphism in Proposition
\ref{prop:K_0} preserves the natural order and scaling.
\end{remark}

\section{$C^*$-algebras of ultragraphs}\label{C*alg-sec}

\begin{definition}
Let $\G = (G^0, \G^1, r,s)$ be an ultragraph. A vertex $v\in G^0$ is
said to be {\em regular} if $0 < |s^{-1}(v) | < \infty$. The set of
all regular vertices is denoted by $G^0_{\rg}\subset G^0$.
\end{definition}

\begin{definition} \label{dfn:CK-G-fam}
For an ultragraph $\G = (G^0, \G^1, r,s)$, a \emph{Cuntz-Krieger
$\G$-family} is a representation $\{p_A\}_{A\in \G^0}$ of $\G^0$ in a
\Ca $\B$ and a collection of partial isometries $\{s_e\}_{e \in \G^1}$ in
$\B$ with mutually orthogonal ranges that satisfy
\begin{enumerate}
\item $s_e^*s_e = p_{r(e)}$ for all $e \in \G^1$,
\item $s_es_e^* \leq p_{s(e)}$ for all $e \in \G^1$,
\item $p_v = \sum_{s(e) = v} s_es_e^*$ for all $v\in G^0_{\rg}$,
\end{enumerate}
where we write $p_v$ in place of $p_{ \{ v \} }$ for $v\in G^0$.

The \Ca $C^*(\G)$ is the \Ca generated by a universal Cuntz-Krieger
$\G$-family $\{p_A, s_e\}$.
\end{definition}

We will show that this definition of $C^*(\G)$ and the following
natural generalization of the definition of Exel-Laca algebras in
\cite{EL} are both equivalent to the original definition of $C^*(\G)$
in \cite[Definition 2.7]{Tom}.

\begin{definition}\label{dfn:EL-G-fam}
For an ultragraph $\G = (G^0, \G^1, r,s)$, an \emph{Exel-Laca
$\G$-family} is a collection of projections $\{p_v\}_{v\in G^0}$ and
partial isometries $\{s_e\}_{e \in \G^1}$ with mutually orthogonal
ranges for which
\begin{enumerate}
\item the collection $\{p_v\}_{v\in G^0} \cup \{s_e^*s_e\}_{e\in \G^1}$
satisfies Condition~(EL),
\item $s_es_e^* \leq p_{s(e)}$ for all $e \in \G^1$,
\item $p_v = \sum_{s(e) = v} s_es_e^*$ for $v\in G^0_{\rg}$.
\end{enumerate}
\end{definition}

\begin{proposition}\label{prop:CK=EL}
For each Cuntz-Krieger $\G$-family $\{p_A, s_e\}$, the collection
$\{p_v,s_e\}$ is an Exel-Laca $\G$-family. Conversely, for each
Exel-Laca $\G$-family $\{p_v,s_e\}$, there exists a unique
representation $\{p_A\}$ of $\G^0$ on the \Ca generated by
$\{p_v,s_e\}$ such that $p_{\{v\}}=p_v$ for $v\in G^0$ and $\{p_A,
s_e\}$ is a Cuntz-Krieger $\G$-family.
\end{proposition}

\begin{proof}
This follows from Proposition \ref{prop:rep<->EL}.
\end{proof}

\begin{corollary}\label{cor:univEL}
Let $\{p_v,s_e\}$ be the Exel-Laca $\G$-family in $C^*(\G)$. For an
Exel-Laca $\G$-family $\{P_v,S_e\}$ on a \Ca $\B$, there exists a
\shom $\phi\colon C^*(\G) \to\B$ such that $\phi(p_v)=P_v$ and
$\phi(s_e)=S_e$. The \shom $\phi$ is injective if $P_v\neq 0$ for all
$v\in G^0$ and there exists a strongly continuous action $\beta$ of
$\T$ on $\B$ such that $\beta_z(P_v)=P_v$ and $\beta_z(S_e)=zS_e$ for
$v\in G^0$, $e\in \G^1$, and $z \in \T$.
\end{corollary}

\begin{proof}
The first part follows from Proposition \ref{prop:CK=EL}, and the
latter follows from \cite[Theorem 6.8]{Tom} because $\phi(p_A) \neq
0$ for all non-empty $A$ if $\phi(p_v)=P_v \neq 0$ for all $v\in
G^0$.
\end{proof}

It is easy to see that Proposition~\ref{prop:CK=EL} is still true if
we replace $\G^0$ by any lattice contained in $\G^0$ and containing
$\{v\}$ and $r(e)$ for all $v \in G^0$ and $e \in \G^1$ (see Lemma
\ref{lem:G'->(*)}). Hence the restriction gives a natural bijection
from Cuntz-Krieger $\G$-families in the sense of Definition
\ref{dfn:CK-G-fam} to the Cuntz-Krieger $\G$-families of
\cite[Definition 2.7]{Tom}. Thus the \Ca $C^*(\G)$ is naturally
isomorphic to the \Ca defined in \cite[Theorem 2.11]{Tom}.

\section{The spectrum of the
commutative $C^*$-algebra $\A_\G$} \label{OmegaG-sec}

Let $\G = (G^0, \G^1, r, s)$ be an ultragraph. In this section, we
describe the spectrum of the commutative \Ca $\A_\G$ concretely. Fix
a listing $\G^1 = \{ e_i \}^\infty_{i=1}$ of the edges of $\G$. As
described in the paragraph following the proof of
Lemma~\ref{lem:G'->(*)}, the \Ca $\A_\G$ is equal to the
inductive limit of the increasing family
$\{\A_\G^{(n)}\}_{n=1}^\infty$, where $\A_\G^{(n)}$ is the \Csa of
$\A_\G$ generated by $\{\delta_{v}:v\in G^0\} \cup
\{\chi_{r(e_i)}:i=1,2,\ldots,n\}$. In order to compute the spectrum
of $\A_\G$, we first compute the spectrum of $\A_\G^{(n)}$ for a
positive integer $n$.

\begin{definition}\label{def:OGn}
For $n \in \N\setminus\{0\}$, we define a topological space
$\Omega_\G^{(n)}$ such that $\Omega_\G^{(n)}=G^0 \sqcup \Delta_n$ as
a set and $A\sqcup Y$ is open in $\Omega_\G^{(n)}$ for $A\subset G^0$
and $Y\subset \Delta_n$ if and only if $|r(\omega)\setminus
A|<\infty$ for all $\omega\in Y$.
\end{definition}

For each $v\in G^0$, $\{v\}$ is open in $\Omega_\G^{(n)}$, and a
fundamental system of neighborhoods of $\omega\in\Delta_n\subset
\Omega_\G^{(n)}$ is
\[
\{A\sqcup\{\omega\} : A\subset G^0, |r(\omega)\setminus A|<\infty\}.
\]
Hence $G^0$ is a discrete dense subset of $\Omega_\G^{(n)}$. Note
that $\Omega_\G^{(n)}$ is a disjoint union of the finitely many
compact open subsets $r(\omega)\sqcup\{\omega\}$ for $\omega\in
\Delta_n$ and the discrete set $G^0\setminus \bigcup_{\omega\in
\Delta_n}r(\omega)$. This fact and the paragraph following
Proposition~\ref{prop:A-G-n generated} show the following:

\begin{lemma}\label{lem:AGn->C}
There exists an isomorphism $\pi^{(n)}\colon \A_\G^{(n)}\to
C_0(\Omega_\G^{(n)})$ such that $\pi^{(n)}(\delta_v)=\delta_v$ and
$\pi^{(n)}(\chi_{r(\omega)})=\chi_{r(\omega)\sqcup\{\omega\}}$ for
$v\in G^0$ and $\omega\in \Delta_n$.
\end{lemma}

\begin{lemma}\label{lem:closure}
For $i=1,2,\ldots,n$, the closure $\overline{r(e_i)}$ of
$r(e_i)\subset \Omega_\G^{(n)}$ is the compact open set
$r(e_i)\sqcup\{\omega\in \Delta_n : \omega_i=1\}$, and we have
$\pi^{(n)}(\chi_{r(e_i)})=\chi_{\overline{r(e_i)}}$.
\end{lemma}

\begin{proof}
This follows from Lemma~\ref{lem:r(omega)}.
\end{proof}

Let $\tDelta_n:=\Delta_n\cup\{0^n\}$. We can define a topology on
$\tOmega_\G^{(n)}:=G^0 \sqcup \tDelta_n$ similarly as in
Definition~\ref{def:OGn} so that $\tOmega_\G^{(n)}$ is the one-point
compactification of $\Omega_\G^{(n)}$. The restriction map
$\{0,1\}^{n+1}\to \{0,1\}^n$ induces a map $\tDelta_{n+1}\to
\tDelta_n$, and hence a map $\tOmega_\G^{(n+1)}\to \tOmega_\G^{(n)}$.
It is routine to check that this map is a continuous surjection, and
the induced \shom $C_0(\Omega_\G^{(n)})\to C_0(\Omega_\G^{(n+1)})$
coincides with the inclusion $\A_\G^{(n)}\hookrightarrow
\A_\G^{(n+1)}$ via the isomorphisms in Lemma~\ref{lem:AGn->C}.

For each element
\[
\omega=(\omega_1,\omega_2,\ldots,\omega_i,\ldots)\in\{0,1\}^\infty,
\]
we define $\omega|_n\in \{0,1\}^{n}$ by
$\omega|_n=(\omega_1,\omega_2,\ldots,\omega_n)$. The space
$\{0,1\}^\infty$ is a compact space with the product topology, and it
is homeomorphic to $\varprojlim \{0,1\}^{n}$.

\begin{definition}
We define
\[
\tDelta_\infty :=\big\{\omega \in\{0,1\}^\infty: \text{$\omega|_n \in
\tDelta_n$ for all $n$}\big\},
\]
and $\Delta_\infty:=\tDelta_\infty\setminus\{0^\infty\}$ where
$0^\infty:=(0,0,\ldots)\in \{0,1\}^\infty$.
\end{definition}

Since $\tDelta_\infty$ is a closed subset of $\{0,1\}^\infty$, the
space $\Delta_\infty$ is locally compact, and its one-point
compactification is homeomorphic to $\tDelta_\infty$. By definition,
$\tDelta_\infty\cong \varprojlim\tDelta_n$.

\begin{definition}
We define a topological space $\Omega_\G$ as follows: $\Omega_\G=G^0
\sqcup \Delta_\infty$ as a set and $A\sqcup Y$ is open in $\Omega_\G$
for $A\subset G^0$ and $Y\subset \Delta_\infty$ if and only if for
each $\omega\in Y$ there exists an integer $n$ satisfying
\begin{enumerate}\rom
\item if $\omega'\in \Delta_\infty$ and $\omega'|_n
= \omega|_n$, then $\omega' \in Y$; and
\item $|r(\omega|_{n}) \setminus A|<\infty$.
\end{enumerate}
\end{definition}

Equivalently, $A\sqcup Y\subset \Omega_\G$ is closed if and only if
$Y\subset \Delta_\infty$ is closed in the product topology on
$\{0,1\}^{\infty}$, and for each $\omega \in \Delta_\infty$ such that
$|r(\omega|_{n})\cap A|=\infty$ for all $n$, we have $\omega \in Y$.

We can define a topology on $\tOmega_\G:=G^0 \sqcup \tDelta_\infty$
similarly as in the definition above, so that $\tOmega_\G$ is the
one-point compactification of $\Omega_\G$ and $\tOmega_\G\cong
\varprojlim\tOmega_\G^{(n)}$.

\begin{lemma}\label{lem:compactopen}
In the space $\Omega_\G$, the closure $\overline{r(e_i)}$ of
$r(e_i)\subset \Omega_\G$ is the compact open set
$r(e_i)\sqcup\{\omega\in \Delta_\infty : \omega_i=1\}$.
\end{lemma}
\begin{proof}
This follows from the homeomorphism $\tOmega_\G\cong
\varprojlim\tOmega_\G^{(n)}$ combined with Lemma~\ref{lem:closure}.
\end{proof}

\begin{proposition}\label{prop:AG->C}
There exists an isomorphism $\pi\colon \A_\G\to C_0(\Omega_\G)$ such
that $\pi(\delta_v)=\delta_v$ and
$\pi(\chi_{r(e)})=\chi_{\overline{r(e)}}$ for $v\in G^0$ and $e\in
\G^1$.
\end{proposition}
\begin{proof}
Taking the inductive limit of the isomorphisms $\pi^{(n)}$ in
Lemma~\ref{lem:AGn->C} produces an isomorphism
\[\textstyle
\pi\colon \A_\G\to \varinjlim C_0(\Omega_\G^{(n)}) \cong
C_0(\Omega_\G).
\]
This isomorphism satisfies the desired condition by
Lemma~\ref{lem:compactopen}.
\end{proof}

By the isomorphism $\pi$ in the proposition above, we can identify
the spectrum of $\A_\G$ with the space $\Omega_\G$.

\section{Topological quivers and $K$-groups}
\label{TopQu&K-sec}

In this section we will construct a topological quiver $\cQ(\G)$ from
$\G$, and show that the \Ca $C^*(\G)$ is isomorphic to the \Ca
$C^*(\cQ(\G))$ of \cite{MT}. Fix an ultragraph $\G = (G^0, \G^1,
r,s)$, and define
\[
    \cQ(\G) := (E(\G)^0, E(\G)^1, r_\cQ, s_\cQ, \lambda_\cQ)
\]
as follows. Let $E(\G)^0:=\Omega_\G$ and
\[
E(\G)^1 := \{(e, x) \in \G^1 \times \Omega_\G : x \in
\overline{r(e)}\},
\]
where $\G^1$ is considered as a discrete set, and
$\overline{r(e)}\subset \Omega_\G$ are compact open sets (see Lemma
\ref{lem:compactopen}).

We define a local homeomorphism $r_\cQ \colon E(\G)^1\to E(\G)^0$ by
$r_\cQ(e,x) := x$, and a continuous map $s_\cQ \colon E(\G)^1\to
E(\G)^0$ by $s_\cQ(e,x) := s(e)\in G^0\subset E(\G)^0$. Since $r_\cQ$
is a local homeomorphism, we have that $r_\cQ^{-1}(x)$ is discrete
and countable for each $x\in E(\G)^0$. For each $x \in E(\G)^0$ we
define the measure $\lambda_x$ on $r_\cQ^{-1}(x)$ to be counting
measure, and set $\lambda_\cQ = \{\lambda_x : x \in E(\G)^0\}$.

Reversing the roles of the range and source maps, we can also regard
$\cQ(\G)$ as a topological graph $E(\G)$ in the sense of \cite{Kat},
and its \Ca $\cO(E(\G))$ is naturally isomorphic to $C^*(\cQ(\G))$
(see \cite[Example~3.19]{MT}). Since some of the results about
$C^*(\cQ(\G))$ which we wish to apply have only been proved in the
setting of \cite{Kat} to date, we will frequently reference these
results; the reversal of edge-direction involved in regarding
$\cQ(\G)$ as a topological graph is implicit in these statements. We
have opted to use the notation and conventions of \cite{MT}
throughout, and to reference the results of \cite{MT} where possible
only because the edge-direction conventions there agree with those
for ultragraphs \cite{Tom}.

We let $E(\G)^0_{\rg}$ denote the largest open subset of $E(\G)^0$
such that the restriction of $s_\cQ$ to $s_\cQ^{-1}(E(\G)^0_{\rg})$
is surjective and proper.

\begin{lemma}\label{E0rg}
We have $E(\G)^0_{\rg}=G^0_{\rg}$.
\end{lemma}

\begin{proof}
Since the image of $s_\cQ$ is contained in $G^0\subset E(\G)^0$, we
have $E(\G)^0_{\rg}\subset G^0$. For each $v\in G^0$, we see that
$v\in E(\G)^0_{\rg}$ if and only if $s_\cQ^{-1}(v)$ is non-empty and
compact because $\{v\}$ is open in $E(\G)^0$. Since $\{e\}\times
\overline{r(e)}\subset E(\G)^1$ is compact for all $e\in\G^1$, it
follows that $s_\cQ^{-1}(v) = \{(e, x)\in E(\G)^1: s(e) = v\}$ is
non-empty and compact if and only if $\{e \in \G^1: s(e) = v\}$ is
non-empty and finite; that is, if and only if $v\in G^0_{\rg}$.
\end{proof}

For the statement of the following theorem, we identify $C_0(E(\G)^0)$
with $\A_\G$ via the isomorphism in Proposition~\ref{prop:AG->C},
and denote by $\chi_e\in C_c(E(\G)^1)$ the characteristic function of the
compact open subset $\{e\}\times \overline{r(e)}\subset E(\G)^1$ for
each $e\in\G^1$. We denote by $(\psi_{\cQ(\G)}, \pi_{\cQ(\G)})$ the
universal generating $\cQ(\G)$-pair, and by $\{p^\G_A, s^\G_e : A \in \G^0,
e \in \G^1\}$ the universal generating Cuntz-Krieger $\G$-family.

\begin{theorem}\label{ultra=topgrah}
There is an isomorphism from $C^*(\G)$ to $C^*(\cQ(\G))$ which is canonical in the sense that it
takes $p^\G_A$ to $\pi_{\cQ(\G)}(\chi_A)$ and $s^\G_e$ to
$\psi_{\cQ(\G)}(\chi_e)$ for all for all $A \in \G^0$ and $e \in \G^1$.
Moreover, this isomorphism is equivariant for the gauge actions on $C^*(\G)$ and
$C^*(\cQ(\G))$.
\end{theorem}
\begin{proof}
It is easy to check using Lemma~\ref{E0rg} and
Proposition~\ref{prop:CK=EL} that:
\begin{enumerate}
\item for each Cuntz-Krieger $\G$-family
$\{p_A, s_e : A \in \G^0, e \in \G^1\}$ there is a unique
$\cQ(\G)$-pair $(\pi_{q,t}, \psi_{p,s})$ satisfying
$\pi_{p,s}(\chi_A) = p_A$ for each $A \in \G^0$ and
$\psi_{p,s}(\chi_e) = s_e$ for each $e \in \G^1$; and
\item for each $\cQ(\G)$-pair $(\pi, \psi)$, the formulae
$p^{\pi,\psi}_A := \pi(\chi_A)$ and $s^{\pi,\psi}_e := \psi(\chi_e)$
determine a Cuntz-Krieger $\G$-family $\{p^{\pi, \psi}_A,
s^{\pi,\psi}_e : A \in \G^0, e \in \G^1\}$.
\end{enumerate}
The result then follows from the universal properties of the two
\CA s $C^*(\G)$ and $C^*(\cQ(\G))$.
\end{proof}

\begin{remark}
To prove Theorem~\ref{ultra=topgrah}, one could alternatively use the
gauge-invariant uniqueness theorems of ultragraphs
\cite[Theorem~6.8]{Tom} or of topological graphs
\cite[Theorem~4.5]{Kat}, or of topological quivers
\cite[Theorem~6.16]{MT}.
\end{remark}

\begin{theorem}
Let $\G = (G^0, \G^1, r, s)$ be an ultragraph. We define
$\partial\colon \Z^{G^0_{\rg}}\to Z_\G$ by
$\partial(\delta_{v})=\delta_{v}-\sum_{e\in s^{-1}(v)}\chi_{r(e)}$
for $v\in G^0_{\rg}$. Then we have $K_0(C^*(\G))\cong
\coker(\partial)$ and $K_1(C^*(\G))\cong \ker(\partial)$.
\end{theorem}

\begin{proof}
Since $C_0(E(\G)^0_{\rg})\cong c_0(G^0_{\rg})$ and $C_0(E(\G)^0)\cong
\A_\G$, we have $K_0(C_0(E(\G)^0_{\rg}))\cong\Z^{G^0_{\rg}}$,
$K_0(C_0(E(\G)^0))\cong Z_\G$ and
$K_1(C_0(E(\G)^0_{\rg}))=K_1(C_0(E(\G)^0))=0$ by Proposition
\ref{prop:K_0}. It is straightforward to see that the map
$[\pi_r]\colon K_0(C_0(E(\G)^0_{\rg}))\to K_0(C_0(E(\G)^0))$ in
\cite[Corollary~6.10]{Kat} satisfies $[\pi_r](\delta_{v})=\sum_{e\in
s^{-1}(v)}\chi_{r(e)}$ for $v\in G^0_{\rg}$. Hence the conclusion
follows from \cite[Corollary 6.10]{Kat}.
\end{proof}

\section{Gauge-invariant ideals} \label{Ideal-sec}

In this section we characterize the gauge-invariant ideals of
$C^*(\G)$ for an ultragraph $\G$ in terms of combinatorial data
associated to $\G$.

For the details of the following, see \cite{MT}. Let $\cQ = (E^0,
E^1, r, s, \lambda)$ be a topological quiver. We say that a subset $U
\subset E^0$ is \emph{hereditary} if, whenever $e \in E^1$ satisfies
$s(e) \in U$, we have $r(e) \in U$. We say that $U$ is
\emph{saturated} if, whenever $v \in E^0_\rg$ satisfies $r(s^{-1}(v))
\subset U$, we have $v \in U$.

Suppose that $U \subset E^0$ is open and hereditary. Then
\[
\cQ_U := (E^0_U, E^1_U, r|_{E^1_U}, s|_{E^1_U},\lambda|_{E^0_U})
\]
is a topological quiver, where $E^0_U=E^0 \setminus U$ and $E^1_U=E^1
\setminus r^{-1}(U)$.

We say that a pair $(U,V)$ of subsets of $E^0$ is \emph{admissible}
if
\begin{enumerate}
\item $U$ is a saturated hereditary open subset of $E^0$; and
\item $V$ is an open subset of $E^0_U$ with
$E^0_\rg \setminus U \subset V \subset (E^0_U)_\rg$.
\end{enumerate}
It follows from \cite[Theorem~8.22]{MT} that the gauge-invariant
ideals of $C^*(\cQ)$ are in bijective correspondence with the
admissible pairs $(U,V)$ of $\cQ$.

Let $\G = (G^0, \G^1, r, s)$ be an ultragraph. We define admissible
pairs of $\G$ in a similar way as above, and show that these are in
bijective correspondence with the gauge-invariant ideals of
$C^*(\G)$.

\begin{definition} \label{def:ideal with new def}
A subcollection $\cH \subset \G^0$ is said to be an {\em ideal} if it
satisfies
\begin{enumerate}
\item $A,B\in\cH$ implies $A\cup B\in\cH$,
\item $A\in\cH$, $B\in\G^0$ and $B\subset A$ imply $B\in\cH$.
\end{enumerate}
\end{definition}

Let $\pi\colon \A_\G\to C_0(\Omega_\G)$ be the isomorphism in
Proposition~\ref{prop:AG->C}. For an ideal $\cH$ of $\G^0$, the set
$\cspa\{\chi_A : A\in \cH\}$ is an ideal of the \Ca $\A_\G$. Hence
there exists an open subset $U_\cH$ of $\Omega_\G$ with
\[
C_0(U_\cH)=\pi\big(\cspa\{\chi_A : A\in \cH\}\big).
\]

\begin{lemma}\label{lem:ideal corresp}
The correspondence $\cH\mapsto U_\cH$ is a bijection from the set of
all ideals of $\G^0$ to the set of all open subsets of $\Omega_\G$.
\end{lemma}

\begin{proof}
Since $\A_\G$ is an AF-algebra, every ideal of $\A_\G$ is generated
by its projections. From this fact, we see that $\cH\mapsto
\cspa\{\chi_A : A\in \cH\}$ is a bijection from the set of all ideals
of $\G^0$ to the set of all ideals of $\A_\G$. Hence the conclusion
follows from the well-known fact that $U\mapsto C_0(U)$ is a
bijection from the set of all open subsets of $\Omega_\G$ to the set
of all ideals of $C_0(\Omega_\G)$.
\end{proof}

\begin{remark}
The existence of this bijection is one of the advantages of changing
the definition of $\G^0$ from that given in \cite{Tom}.
\end{remark}

\begin{lemma}\label{lem:v and r(e)}
Let $\cH$ be an ideal of $\G^0$, and let $U_\cH\subset \Omega_\G$ be
the corresponding open set. Then for $v\in G^0$, $\{v\}\in \cH$ if
and only if $v\in U_\cH$, and for $e\in\G^1$, $r(e)\in \cH$ if and
only if $\overline{r(e)}\subset U_\cH$.
\end{lemma}
\begin{proof}
This follows from Proposition~\ref{prop:AG->C}.
\end{proof}

\begin{definition}\label{dfn:sat,hered}
We say that an ideal $\cH \subset \G^0$ is \emph{hereditary} if,
whenever $e \in \G^1$ satisfies $\{s(e)\} \in \cH$, we have $r(e) \in
\cH$, and that it is \emph{saturated} if, whenever $v \in G^0_\rg$
satisfies $r(e) \in \cH$ for all $e \in s^{-1}(v)$, we have $\{v\}
\in \cH$.
\end{definition}

\begin{proposition}\label{prop:H<->U_H}
An ideal $\cH$ of $\G^0$ is hereditary (resp.\ saturated) if and only
if the corresponding open subset $U_\cH\subset \Omega_\G=E(\G)^0$ is
hereditary (resp.\ saturated) in the topological quiver $\cQ(\G)$.
\end{proposition}
\begin{proof}
An open subset $U\subset \Omega_\G=E(\G)^0$ is hereditary if and only
if, whenever $(e,x)\in E(\G)^1$ satisfies $s_\cQ(e,x) = s(e)\in U$,
we have $r_\cQ(e,x) = x\in U$. This is equivalent to the statement
that, whenever $e\in\G^1$ satisfies $s(e)\in U$, we have
$\overline{r(e)}\subset U$. Thus Lemma~\ref{lem:v and r(e)} shows
that an ideal $\cH$ is hereditary if and only if $U_\cH$ is
hereditary.

An open subset $U\subset \Omega_\G=E(\G)^0$ is saturated if and only
if, whenever $v\in E(\G)^0_{\rg}=G^0_{\rg}$ satisfies
$r_\cQ(s_\cQ^{-1}(v))\subset U$, we have $v\in U$. For $v\in
G^0_{\rg}$, we have $r_\cQ(s_\cQ^{-1}(v))=\bigcup_{e \in
s^{-1}(v)}\overline{r(e)}$. Hence $U$ is saturated if and only if,
whenever $v\in G^0_{\rg}$ satisfies $\overline{r(e)}\subset U$ for
all $e \in s^{-1}(v)$, we have $v\in U$. Thus Lemma~\ref{lem:v and
r(e)} again shows that an ideal $\cH$ is saturated if and only if
$U_\cH$ is saturated.
\end{proof}

\begin{definition}\label{def:G0Hrg}
Let $\cH$ be a hereditary ideal of $\G^0$. For $v\in G^0$, we define
$s_{\G/\cH}^{-1}(v)\subset \G^1$ by
\[
s_{\G/\cH}^{-1}(v):= \{e\in\G^1:\text{$s(e)=v$ and
$r(e)\notin\cH$}\}.
\]
We define $(G^0_{\cH})_\rg\subset G^0$ by
\[
(G^0_{\cH})_\rg:= \{v\in G^0: \text{$s_{\G/\cH}^{-1}(v)$ is non-empty
and finite}\}.
\]
\end{definition}

Since $\cH$ is hereditary, if $\{v\}\in \cH$ then we have
$s_{\G/\cH}^{-1}(v)=\emptyset$ and hence $v\notin (G^0_{\cH})_\rg$.

\begin{lemma}
A hereditary ideal $\cH$ of $\G^0$ is saturated if and only if,
whenever $v\in G^0_\rg$ satisfies $\{v\}\notin \cH$, we have $v\in
(G^0_{\cH})_\rg$.
\end{lemma}
\begin{proof}
An element $v\in G^0_\rg$ is in $(G^0_{\cH})_\rg$ if and only if
$s_{\G/\cH}^{-1}(v)\subset s^{-1}(v)$ is non-empty, which occurs if
and only if there exists $e\in s^{-1}(v)$ with $r(e)\notin \cH$.
\end{proof}

Let $\cH$ be a hereditary ideal of $\G^0$, and $U_\cH\subset
\Omega_\G=E(\G)^0$ be the corresponding open subset which is
hereditary by Proposition~\ref{prop:H<->U_H}. As in the beginning of
this section, we obtain a topological quiver $\cQ(\G)_{U_\cH}$.

\begin{lemma}\label{lem:V sets}
We have $(E(\G)^0_{U_\cH})_\rg=(G^0_{\cH})_\rg$.
\end{lemma}
\begin{proof}
Since the image of $s_\cQ|_{E(\G)^1_{U_\cH}}$ is contained in $G^0
\setminus U_\cH$, we have $(E(\G)^0_{U_\cH})_\rg\subset G^0 \setminus
U_\cH$. For $v\in G^0$, $v\in U_\cH$ implies $\{v\}\in H$ by
Lemma~\ref{lem:v and r(e)}, and this implies $v\notin
(G^0_{\cH})_\rg$ as remarked after Definition~\ref{def:G0Hrg}. Hence
we have $(G^0_{\cH})_\rg\subset G^0 \setminus U_\cH$. An element
$v\in G^0 \setminus U_\cH$ is in $(E(\G)^0_{U_\cH})_\rg$ if and only
if $s_\cQ^{-1}(v)\cap E(\G)^1_{U_\cH}$ is non-empty and compact
because $\{v\}$ is open in $E(\G)^0_{U_\cH}$. Since
\[
s_\cQ^{-1}(v)\cap E(\G)^1_{U_\cH} =\{(e, x)\in E(\G)^1: \text{$s(e) =
v$ and $x\notin U_\cH$}\},
\]
$s_\cQ^{-1}(v)\cap E(\G)^1_{U_\cH}$ is non-empty and compact if and
only if
\[
\{e\in\G^1: \text{$s(e) = v$ and $\overline{r(e)}\not\subset
U_\cH$}\}
\]
is non-empty and finite. This set is equal to $s_{\G/\cH}^{-1}(v)$ by
Lemma~\ref{lem:v and r(e)}. Therefore an element $v\in G^0 \setminus
U_\cH$ is in $(E(\G)^0_{U_\cH})_\rg$ if and only if $v\in
(G^0_{\cH})_\rg$. Thus $(E(\G)^0_{U_\cH})_\rg=(G^0_{\cH})_\rg$ as
required.
\end{proof}

By Lemma~\ref{lem:V sets}, the subset $(E(\G)^0_{U_\cH})_\rg\subset
E(\G)^0_{U_\cH}$ is discrete.

\begin{definition}
Let $\G = \{G^0, \G^1, r, s\}$ be an ultragraph. We say that a pair
$(\cH,V)$ consisting of an ideal $\cH$ of $\G^0$ and a subset $V$ of
$G^0$ is \emph{admissible} if $\cH$ is hereditary and saturated and
$V\subset (G^0_{\cH})_\rg\setminus G^0_\rg$.
\end{definition}

\begin{definition}
For an admissible pair $(\cH,V)$ of an ultragraph $\G$, we define an
ideal $I_{(\cH, V)}$ of $\cC^*(\G)$ to be the ideal generated by the
projections
\[
\{p_A : A \in \cH\}\cup \Big\{p_v - \sum_{e \in s_{\G/\cH}^{-1}(v)}
s_e s^*_e : v\in V\Big\}.
\]
For an ideal $I$ of $\cC^*(\G)$, we define $\cH_I:=\{A \in \G^0 : p_A
\in I\}$ and
\[
V_I:=\Big\{v\in (G^0_{\cH_I})_\rg\setminus G^0_\rg : p_v-\sum_{e\in
s_{\G/\cH_I}^{-1}(v)} s_e s^*_e\in I\Big\}.
\]
\end{definition}

\begin{theorem}\label{thm:all gi ideals}
Let $\G$ be an ultragraph. Then the correspondence $I\mapsto (\cH_I,
V_I)$ is a bijection from the set of all gauge-invariant ideals of
$C^*(\G)$ to the set of all admissible pairs of $\G$, whose inverse
is given by $(\cH, V)\mapsto I_{(\cH, V)}$.
\end{theorem}
\begin{proof}
By Theorem~\ref{ultra=topgrah}, the gauge-invariant ideals of
$C^*(\G)$ are in bijective correspondence with the gauge-invariant
ideals of $C^*(\cQ(\G))$. We know that the latter are indexed by
admissible pairs $(U,V)$ of $\cQ(\G)$ by \cite[Theorem~8.22]{MT}.
Proposition~\ref{prop:H<->U_H} and Lemma~\ref{lem:V sets} show that
$(\cH, V)\mapsto (U_\cH,V\cup (G^0_\rg\setminus U_\cH))$ is a
bijection from the set of all admissible pairs of $\G$ to the one of
$\cQ(\G)$. Thus we get bijective correspondences between the set of
all gauge-invariant ideals of $C^*(\G)$ and the set of all admissible
pairs of $\G$. By keeping track of the arguments in
\cite[Section~8]{MT}, we see that the bijective correspondences are
given by $I\mapsto (\cH_I, V_I)$ and $(\cH, V)\mapsto I_{(\cH, V)}$.
\end{proof}

\begin{remark}
The theorem above naturally generalizes \cite[Theorem~3.6]{BHRS}.
\end{remark}

\section{Condition (K)} \label{CondK-sec}

In this section we define a version of Condition~(K) for ultragraphs,
and show that this condition characterizes ultragraphs $\G$ such that
every ideal of $C^*(\G)$ is gauge-invariant.  

Let $\G$ be an ultragraph. For $v \in G^0$, a \emph{first-return path
based at $v$} in $\G$ is a path $\alpha = e_1 e_2 \cdots e_n$ such
that $s(\alpha) = v$, $v \in r(\alpha)$, and $s(e_i)\neq v$ for $i =
2,3,\ldots,n$. When $\alpha$ is a first-return path based at $v$, we
say that $v$ \emph{hosts the first-return path $\alpha$}.

Note that there is a subtlety here: a first-return path based at $v$
may pass through other vertices $w \not= v$ more than once (that is,
we may have $s(e_i) = s(e_j)$ for some $1 < i,j \leq n$ with $i \not=
j$), but no edge other than $e_1$ may have source $v$.

\begin{definition}
Let $\G$ be an ultragraph. We say that $\G$ satisfies Condition~(K)
if every $v \in \G^0$ which hosts a first-return path hosts at least
two distinct first-return paths.
\end{definition}

\begin{example}
The graph
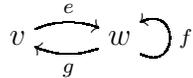
\begin{equation*}
\xymatrix{  v \ar@/^/[r]^e & w \ar@(dr,ur)[]_f  \ar@/^/[l]^g \\}
\end{equation*}
satisfies Condition~(K) because $v$ hosts infinitely many
first-return paths $eg,efg,effg, \ldots$, and $w$ hosts two
first-return paths $f$ and $ge$. Note that all first-return paths
based at $v$ except $eg$ pass through the vertex $w$ more than once.
\end{example}

\begin{proposition}
Let $\G = (G^0, \G^1, r, s)$ be an ultragraph. Then every ideal of
$C^*(\G)$ is gauge-invariant if and only if $\G$ satisfies
Condition~(K).
\end{proposition}
\begin{proof}
In the same way as above, we can define first-return paths in
the topological graph $E(\G)$. It is straightforward to see
that for each $v\in G^0$, first-return paths $\alpha = e_1 e_2
\cdots e_n$ based at $v$ in $\G$ correspond bijectively to
first-return paths
\[
l=(e_1, s(e_2))(e_2, s(e_3)) \cdots (e_n, s(e_1))
\]
based at $v\in G^0\subset E(\G)^0$ in $E(\G)$.

Recall (see \cite[Definition~7.1]{Kat2} and the subsequent
paragraph for details) that $\operatorname{Per}(E(\G))$ denotes
the collection of vertices $v \in E(\G)^0$
such that $v$ hosts exactly one first-return path in $E(\G)$,
and $v$ is isolated in
\[
\{s_{\cQ}(l) : \text{$l$ is a path in $E(\G)$ with $r_{\cQ}(l) = v$}\}
\]
(recall that the directions of paths are reversed
when passing from the quiver $\cQ(\G)$ to the topological graph
$E(\G)$).
Wee see that \cite[Theorem~7.6]{Kat2} implies that every ideal of
$\cO(E(\G))$ is gauge invariant if and only if
$\operatorname{Per}(E(\G))$ is empty. 
Since the isomorphism of $C^*(\G)$ with $\cO(E(\G))$ is gauge equivariant, 
it therefore suffices to show that $\operatorname{Per}(E(\G))$ is empty 
if and only if $\G$ satisfies Condition~(K).

The image of $s_\cQ$ is contained in the discrete set $G^0\subset E(\G)^0$. 
Thus $v \in E(\G)^0$ belongs to
$\operatorname{Per}(E(\G))$ if and only if $v \in G^0\subset
E(\G)^0$ and $v$ hosts exactly one first-return path in $E(\G)$.
By the first paragraph of this proof, 
$v \in G^0\subset E(\G)^0$ hosts exactly one first-return path in $E(\G)$ 
if and only if $v$ hosts exactly one first-return path in $\G$.
Hence $\operatorname{Per}(E(\G))$ is empty 
if and only if $\G$ satisfies Condition~(K).
\end{proof}

\end{document}